\documentclass[12pt]{amsart}
\usepackage{amssymb,amsmath}
\numberwithin{equation}{section}

\def\T{\text}
\def\1#1{\overline{#1}}
\def\2#1{\widetilde{#1}}
\def\3#1{\widehat{#1}}
\def\4#1{\mathbb{#1}}
\def\5#1{\frak{#1}}
\def\6#1{{\mathcal{#1}}}
\def\C{{\4C}}

\def\B{\Bbb B}

\def\phi{\varphi}

\newtheorem{Thm}{Theorem}[section]
\newtheorem{Cor}[Thm]{Corollary}
\newtheorem{Pro}[Thm]{Proposition}
\newtheorem{Lem}[Thm]{Lemma}
\theoremstyle{definition}\newtheorem{Def}[Thm]{Definition}
\theoremstyle{remark}
\newtheorem{Rem}[Thm]{Remark}
\newtheorem{Exa}[Thm]{Example}

\def\Label#1{\label{#1}}
\def\bl{\begin{Lem}}
\def\el{\end{Lem}}
\def\bp{\begin{Pro}}
\def\ep{\end{Pro}}
\def\bt{\begin{Thm}}
\def\et{\end{Thm}}
\def\bc{\begin{Cor}}
\def\ec{\end{Cor}}
\def\bd{\begin{Def}}
\def\ed{\end{Def}}
\def\br{\begin{Rem}}
\def\er{\end{Rem}}
\def\be{\begin{Exa}}
\def\ee{\end{Exa}}
\def\bpf{\begin{proof}}
\def\epf{\end{proof}}
\def\ben{\begin{enumerate}}
\def\een{\end{enumerate}}

\def\1alpha{[\frac1\alpha]}
\def\T{\text}

\def\C{{\Bbb C}}

\numberwithin{equation}{section}
\def\T{\text}

\newtheorem{theorem}{Theorem  }[section]
\newtheorem{definition}[theorem]{Definition }
\newtheorem{lemma}[theorem]{Lemma  }
\newtheorem{proposition}[theorem]{Proposition  }
\newtheorem{corollary}[theorem]{Corollary }
\newtheorem{example}[theorem]{\it Example }

\begin{document}
\title[Holomorphic extension from a convex hypersurface]{Holomorphic extension from a convex hypersurface}
\author[L.~Baracco]
{Luca Baracco}
\address{Dipartimento di Matematica, Universit\`a di Padova, via 
Trieste 63, 35121 Padova, Italy}
\email{baracco@math.unipd.it}
\maketitle
\begin{abstract}
We discuss a general result of holomorphic extension of a real analytic function $f$ defined on the boundary $\partial D$ of a real analytic strictly convex subset $D\subset\subset \C^n$. We show that this follows from the hypothesis of separate holomorphic extension along stationary/extremal discs.
\newline
MSC: 32V10, 32N15, 32D10
\end{abstract}
\section{Stationary discs and holomorphic extension }
The problem of testing analyticity on a domain $D\subset \C^n$ by a family of discs has attracted a great deal of work. The first significant result goes back to Stout \cite{S77} who uses as testing family all the straight lines. Reducing the testing family, Agranovsky and Semenov  \cite{AS71} use the lines which meet an open subset $D'\subset\subset D$. It is classical that the lines which meet a single point $z_o\in D$ do not suffice not even in the case of the sphere $\B^n$. 
Other testin families are considered among others, by \cite{GS},\cite{R},\cite{G88}. In the present paper
for a strictly convex $C^\omega $ domain $D$, we prove that the stationary discs  passing through a point of $\bar D$ is a testing family if the point belongs to the boundary $ \partial D$ and, otherwise, if it is supplemented by another $(2n-2)$-parameter, generic family. In particular this second can be chosen as the family of stationary discs through another point of $D$. This result is also present in a recent paper by Agranovsky \cite{A09}.
We deal with stationary/extremal discs in the sense of Lempert \cite{L81}. We first introduce some terminology. A disc $A$ is the holomorphic image of the standard disc $\Delta$;  $\Bbb PT^*\C^n$ is the cotangent bundle with projectivized fibers, and $\pi$ the projection on the base point; $\Bbb T^*_{\partial D}\C^n$ the projectivized conormal bundle to $\partial D$ in $\C^n$. 
\bd
\Label{d1.1}
A disc $A$ of $D$ is said to be stationary when it is endowed with a meromorphic lift $A^*\subset \Bbb PT^*\C^n$ with a simple pole attached to $ T^*_{\partial D}\C^n$, that is, satisfying $\partial A^*\subset \Bbb T^*_{\partial D}\C^n$.
\ed
We fix a stationary disc $A_o$ of $D$ and, in the $\epsilon$-neighborhood of $A_o$, consider a certain number of $(2n-2)$-parameter families of stationary discs $\{\mathcal V_j\}_{j=1,...k}$ smoothly depending on the parameters. We denote by $\mathcal V_j^*$ the family of lifts of the discs in $\mathcal V_j$ and define
\begin{equation}
\Label{1.1}
M_j:=\underset {A^*\in\mathcal V_j^*}\cup A^*.
\end{equation}
The set $M_j$ is  generically a CR manifold with CR dimension 1 except at the points  of a closed set; we denote by $M_j^{\T{reg}}$ the complement of this set. We assume
\begin{equation}
\Label{1.2}
A^*_o\subset\underset j\cup M_j^{\T{reg}}.
\end{equation}
 Here is our main result.
\bt
\Label{t1.1}
Let $D\subset\subset\C^n$ be a strictly convex domain with $C^\omega$ boundary and $f$ a $C^\omega$ function on $\partial D$. Suppose that $f$ extends holomorphically along each disc $A\in\underset{j=1,...,k}\cup \mathcal V_j$ and that the sets $M_j$ which collect the discs of  $\mathcal V_j$ satisfy \eqref{1.2}. Then $f$ extends holomorphically to $D$.
\et
The proof is given in next section. 
Theorem~\ref{t1.1} states a general principle: $(2n-2)$-families of stationary discs generically suffice. To exhibit explicit families the following criterion is very effective. 
\bt
\Label{t1.2} Let $\mathcal V$ be the discs through a point $z_o\in \bar D$ and $M$ the union of their lifts. 
\begin{itemize}
\item[(i)] If $z_o\in D$, then
$$
A_o^*\setminus \pi^{-1}(z_o)\subset M^{\T{reg}}.
$$
\\
\item[(ii)] If $z_o\in \partial D$, then
$$
A_o^*\subset M^{\T{reg}}.
$$
\end{itemize}
\et
\bpf
(i): We first assume that $D$ coincides with the unit ball $\B^n$. It is classical that the stationary discs are the straight lines. By a biholomorphic transformation of $\B^n$ we can displace $z_o$ at $0$. It is helpful to use the parametrization
\begin{equation*}
\begin{matrix}
\partial\B^n\times(0,1)&\to&M
\\
(z,r)&\mapsto &(zr,[\bar z]),
\end{matrix}
\end{equation*}
where brackets denote projectivized coordinates. For fixed $r>0$, this describes a totally real maximal manifold of $\Bbb P T^*\C^n$; thus $\dim_{CR}M\le 1$. On the other hand, $M$ is foliated by discs and therefore $\dim_{CR}M=1$. 

Instead, for $r=0$, we have $TM|_0=\{0\}\times \Bbb P^{n-1}_\C$; thus any point of $M|_0$ is CR singular since there the CR dimension jumps from $1$ to $n-1$. 

We pass now to a general strictly convex domain $D$. We know from \cite{L81} that there is a mapping $\Psi:\Bbb B^n\to D$ which interchanges $0$ with $z_o$, is $C^\omega$ outside $0$, transforms holomorphically the lines of $\Bbb B^n$ (denoted $A_{\B^n}$) into the stationary discs of $D$ through $z_o$ (denoted $A_D$), and which fixes the tangent directions at the ``centers". Therefore, $\Psi$ lifts in a natural way to a mapping between the manifold $M_{\Bbb B^n}$ (the union of the $A^*_{\Bbb B^n}$'s) to the corresponding manifold $M_D$ (the union of $A^*_D$'s). 
Denote by $\B^n_r$ the ball of radius $r$ and put $D_r:=\Psi(\B^n_r)$; we  know from the theory of Lempert that
$$
A^*_{D_r}=(A^*_D)|_{D_r}.
$$
Since $(A^*_{D_r})|_{\partial D_r}\subset \Bbb PT^*_{\partial D_r}\C^n$, it follows that $M_D\setminus\pi^{-1}(z_o)\subset \underset r\cup \Bbb PT^*_{\partial D_r}\C^n$. Thus, $\Bbb PT^*_{\partial D_r}\C^n$ being maximal totally real for any $r$, we conclude that $M|_{\partial D}$ is a CR manifold except at points of $\pi^{-1}(z_o)$ and that it is CR-diffeomorphic, via $\Psi^*$, to $M_{\B^n}\setminus\pi^{-1}(0)$ .

\noindent
(ii): The proof is the same as in (i) but uses the boundary version of the Riemann-Lempert mapping Theorem as in Chang-Hu-Lee \cite{CL88}.

\epf
We are ready for the following, explicit, result.

\bt
\Label{t1.3}
Let $D$ be strictly convex with $C^\omega$ boundary and let $f\in C^\omega(\partial D)$. Either of the following hypothesis is sufficient for holomorphic extension of $f$ to $D$.
\begin{itemize}
\item[(i)] $f$ extends holomorphically along the stationary discs passing through two points of $D$.
\\
\item[(ii)] $f$ extends along the discs through a boundary point of $\partial D$.
\end{itemize}
\et
\bpf
(i): Let $z_o$ and $w_o$ be the ``centers" of the two systems of discs, let $A_o$ be the disc which connects $z_o$ to  $w_o$,  and denote by $M^{z_o}$ and $M^{w_o}$ the union of the discs through $z_o$ and $w_o$ respectively. 
The lift $A^*_o$ is contained in $(M^{z_o})^{\T{reg}}$ apart from a single point over $z_o$; but this point is contained in $M^{w_o}$. Thus we can apply Theorem~\ref{t1.1}.

\noindent
(ii): If $z_o\in\partial D$, we have directly $A^*_o\subset (M^{z_o})^{\T{reg}}$.

\epf
\br
Note that in (ii) the family of discs $\mathcal V^{w_o}$ is only used to cover the singular point of $M^{z_o}$ over $z_o$; for this purpose, a much more general family than of discs through another point $w_o$ is suitable.
\er

\br
Discs by two points of the ball are also present, as a
testing family, in the recent preprint \cite{A09} by Agranovsky.
\er

\section{Proof of Theorem~\ref{t1.1}}

Before starting the proof, we  have to recall the main results from \cite{L81} which will be on use. 
Stationary discs are stable under reparametrization. In particular, the pole can be displaced at any of their interior points. It is convenient to identify the lift $A^*$ to its image  in the projectivized bundle $\mathbb{P} T^*\C^n$ with coordinates $(z,[\zeta])$. We assume that $D$ is strictly convex and that $\partial D\in C^\omega$. In this situation, a stationary disc and its lift $A^*$ are $C^\omega$ up to $\partial \Delta$. Moreover, one has the following basic result for whose proof we refer to \cite{L81}.
\bp
\Label{p1.1}
For any point $(z,[\zeta])\in \mathbb{P} T^*\C^n|_D$ there is unique, up to reparametrization, the stationary disc whose lift $A^*_{(z,[\zeta])}$ contains $(z,[\zeta])$. Moreover, the correspondence 
\begin{equation}
\Label{1.5}
(z,[\zeta])\mapsto A^*_{(z,[\zeta])},\qquad \mathbb{P} T^*\C^n|_D\to C^\omega(\bar\Delta),
\end{equation}
is a $C^\omega$ diffeomorphism.
\ep
 
We begin now the proof of Theorem~\ref{t1.1} and first remark that at any point of $M_j^{\T{reg}}$, the CR  structure is fully provided by the discs $A^*\in V_j$ by which $M_j$ is foliated. Notice that $M_j$ has a natural ``edge" $E_j :=\underset{A\in\mathcal V_j}\cup\partial A^*$. The function $f$ can be naturally lifted to a function $F$ on $M_j$ by gluing the bunch of separate holomorphic extensions $\{f_{A}\}_{A\in \mathcal V_j}$. This is defined by
\begin{equation*}
F(z,[\zeta])=f_{A_{(z,[\zeta])}}(z),
\end{equation*}
where $A_{(z,[\zeta])}$ is the unique stationary disc of $\mathcal V_j$ whose lift $A^*_{(z,[\zeta])}$ passes through $(z,[\zeta])$. The  crucial point here is that the $A$'s may overlap on $\C^n$ but the $A^*$'s do not in $\mathbb{P} T^*\C^n$. The function $F$ is CR on $M_j^{\T{reg}}$. Moreover, since $f\in C^\omega(\partial D)$ and $T^*_{\partial D}\C^n$ is maximally totally real with complexification $\mathbb{P} T^*\C^n$, then $F$ extends holomorphically to a full neighborhood of $\mathbb{P} T^*_{\partial D}\C^n$ in $\mathbb{P} T^*\C^n$. By propagation of holomorphic extendibility on $M_j^{\T{reg}}$ along the discs $A^*_{(z,[\zeta])}$, $F$ extends holomorphically to a neighborhood of $M_j^{\T{reg}}$. Since $A^*_o\subset \underset j\cup M_j^{\T{reg}}$, then we get the conclusion
\begin{equation}
\Label{2.1}
\T{$F$ is holomorphic in a neighborhood of $A^*_o$.}
\end{equation}
We prove now that \eqref{2.1} implies
\begin{equation}
\Label{2.2}
\T{$F$ is holomorphic in a neighborhood of any other stationary disc $A^*_1$ of $D$.}
\end{equation}
To see this, we suppose $A^*_{(z,[\zeta])}(0)=(z,[\zeta])$ and define a function $G$ by means of Cauchy integral
$$
G(z,[\zeta]):=(2\pi)^{-1}\int_{\partial\Delta}\frac{f\circ A_{(z,[\zeta])}(\tau)}\tau d\tau.
$$
This is defined for any $(z,[\zeta])\in \mathbb{P} T^*\C^n|_D$, is real analytic, and satisfies
$$
G=F\qquad\T{in a neighborhood of $A^*_o$.}
$$
Hence $F$, identified to $G$, extends holomorphically to the full $\mathbb{P} T^*\C^n|_D$. Since $\mathbb{P} T^*\C^n|_D$ is covered by the discs $A^*_{(z,[\zeta])}$ for $(z,[\zeta])\in \mathbb{P} T^*\C^n|_D$ (by Proposition~\ref{p1.1}), since $\partial A^*_{(z,[\zeta])}\subset \mathbb{P} T^*_{\partial D}\C^n$ and since $F$ is bounded over these boundaries, then $F$ is in fact bounded on the whole $\mathbb{P} T^*\C^n|_D$. Therefore it is constant with respect to $[\zeta]$. Thus it is a function of $z$ only, the holomorphic extension of $f$ to $D$.

\hskip14cm $\Box$

\end{document}